# Isoperimetric Pentagonal Tilings

Ping Ngai Chung, Miguel A. Fernandez, Yifei Li, Michael Mara, Frank Morgan, Isamar Rosa Plata, Niralee Shah, Luis Sordo Vieira, Elena Wikner

## Abstract

We identify least-perimeter unit-area tilings of the plane by convex pentagons, namely tilings by Cairo and Prismatic pentagons, find infinitely many, and prove that they minimize perimeter among tilings by convex polygons with at most five sides.

## 1. Introduction

In 2001, Thomas Hales ([H]; see [M1, Chapt. 15]) proved the Honeycomb Conjecture, which says that regular hexagons provide the least-perimeter unit-area way to tile the plane. In this paper, we seek perimeter-minimizing tilings of the plane by unit-area pentagons. The regular pentagon has the least perimeter, but it does not tile the plane. There are many planar tilings by a single irregular pentagon or by many different unit-area pentagons; for some simple examples see Figure 3. Which of them has the least average perimeter per tile? Our main Theorem 3.2 proves that the Cairo and Prismatic tilings of Figure 1 minimize perimeter, assuming that the pentagons are convex. We conjecture that this convexity assumption is unnecessary. The Cairo and Prismatic tiles have identical perimeter $2\sqrt{2+\sqrt{3}} \approx 3.86$, less than a unit square's (4) but more than the regular pentagon's ($\approx 3.81$). There are infinitely many other equally efficient tilings by mixtures of Cairo and Prismatic tiles, as in Figure 2 (see also Remark 2.3 and Figures 4−15).

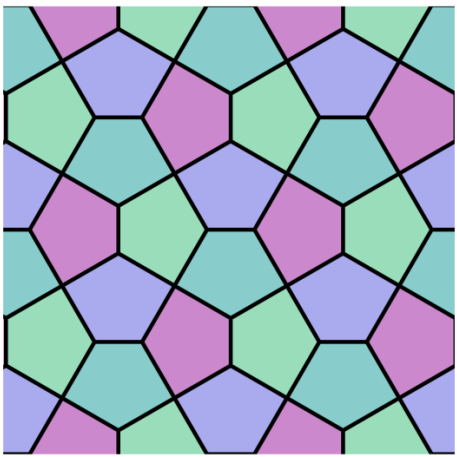
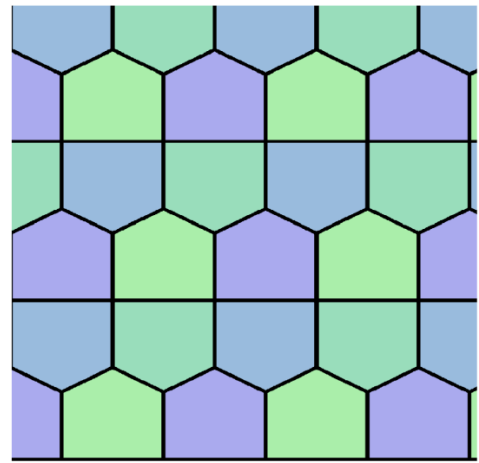

Cairo Pentagonal
Tiling

Prismatic Pentagonal
Tiling

Figure 1. Minimal pentagonal tilings. Each tile has two right angles and three angles of 2π/3 and is circumscribed about a circle.

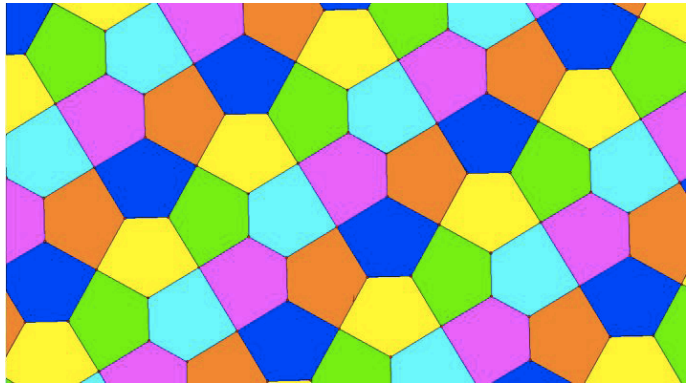

Figure 2. One of infinitely many equally efficient tilings by mixtures of Cairo and Prismatic pentagons.

The tiling of the plane by pentagons has intrigued mathematicians and has been the subject of numerous studies (see the beautiful survey by Schattschneider [S1] and recent work of Sugimoto and Ogawa [SO]). Only in 1985 [HH] were the five types of tilings by a single *equilateral* convex pentagon proved complete. There are 14 known types of tilings by a single convex pentagon, pictured in Figure 3, and more tilings by nonconvex pentagons.

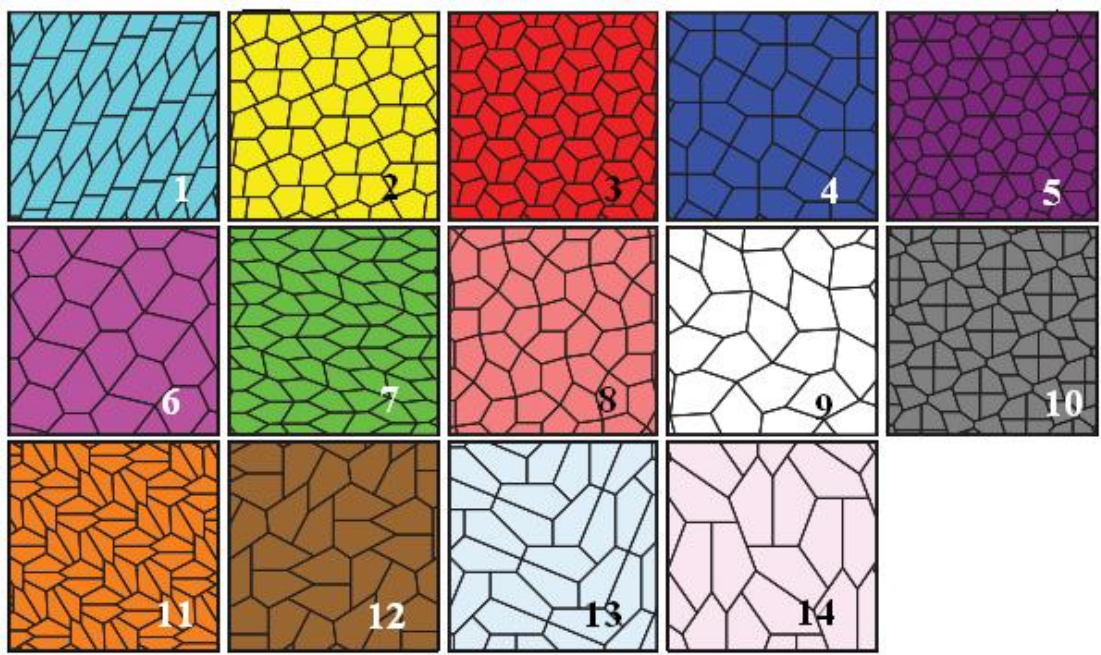

Figure 3. The 14 known types of planar tilings by a convex pentagon.
http://en.wikipedia.org/wiki/Pentagon_tiling

Along the way we prove an isoperimetric result we were not able to find in the literature (see Prop. 3.1 and Fig. 16):

*For n given angles $0 < a_i \leq \pi$ summing to $(n-2)\pi$, the n-gon circumscribed about the unit circle maximizes area for given perimeter.*

This complements the well-known result that for given *edge lengths*, the *n*-gon *inscribed* in a circle maximizes area. (Idea of short proof: otherwise take the inscribed n-gon, with the surmounted circular arcs of the circle, and deform it to enclose more area. With the surmounted circular arcs, it now beats the circle, a contradiction.)

The proof that the Cairo and Prismatic tilings minimize perimeter begins with the Euler characteristic formula, which implies that on average at least three angles per pentagon belong to vertices of the tiling of degree three. Thus at least 3/5 of the angles are "large" angles with average measure $2\pi/3$. Certain convexity arguments show it best for three angles of each pentagon to equal $2\pi/3$ and for the others to equal $\pi/2$, namely, the Cairo and Prismatic tiles. In the problematic case that one or two of the angles is $\pi$, the pentagon becomes a quadrilateral or a triangle and we need to enlarge our category to convex *n*-gons with $n \leq 5$, not necessarily meeting edge to edge. The Euler consequence generalizes to the statement that on the average the number of large angles is at least $3n-12$. Lower bounds for the perimeter contributed separately by triangles, quadrilaterals, and pentagons are compared via linear approximations.

*Acknowledgements.* We thank Doris Schattsneider, Daniel Huson, Max Engelstein, and David Thompson for their help. Morgan thanks the Fields Institute, where he spent fall 2010 on sabbatical. For funding we thank the National Science Foundation for grants to Morgan and to the Williams College "SMALL" REU, Williams College, Berea College, and the Mathematical Association of America for supporting trips to speak at MathFest.

## 2. Tilings

**2.1. Definitions.** A polygonal *tiling* is a decomposition of the plane as a union of polygonal regions which meet only along their boundaries. The tiling is *edge to edge* if the tiles meet only along entire edges. A *monohedral* tiling is a tiling by congruent copies of a single *prototile*.

The *perimeter ratio* of a planar tiling is defined as the limit superior as $R$ approaches infinity of the perimeter inside an $R$-disc about the origin divided by the area $\pi R^2$.

Although the terms are sometimes used more broadly, we define the *Cairo* and *Prismatic* pentagons as unit-area pentagons circumscribed about a circle with two right angles and three angles of $2\pi/3$, as in Figure 1. Proposition 3.1 shows that circumscribed polygons minimize perimeter for given angles. The Cairo tile has one shorter edge of length $a = (2/3)\sqrt{6 - 3\sqrt{3}} \approx .5977$ and four equal longer edges of length $b = (3 + \sqrt{3})\sqrt{2 - \sqrt{3}} / 3 \approx .8165$. The Prismatic tile has two adjacent shorter edges of length $a$ (forming the roof of the house), two edges of length $b$ (the walls), and a still longer base of length $2\sqrt{2 - \sqrt{3}} \approx 1.0353$.

**Proposition 2.2.** *There is a unique planar tiling by Cairo tiles. There is a unique edge-to-edge planar tiling by Prismatic tiles.*

*Proof sketch.* For the Cairo tile, the shorter edges have vertices of degree 3, so every tiling must be edge to edge. Working out from a vertex of degree 4 it is easy to see that an edge-to-edge tiling by the Cairo or the Prismatic tile is unique.

**Remark 2.3 (Nonuniqueness).** The first version of this paper asked whether you can tile the plane with mixtures of Cairo and Prismatic tiles. An MIT undergraduate, shortly recruited as a co-author (Chung), discovered uncountably many tilings by mixtures of Cairo and Prismatic tiles. In Figure 2, there are alternating diagonal rows of Cairo and Prismatic tiles. The Cairo tiles are grouped in hexagons of four, while the Prismatic tiles are grouped in twos. Uncountably many other such tilings may be obtained by placing arbitrarily many copies of rows of one type before copies of rows of the other type. Most of these tilings are nonperiodic. One other example is shown in Figure 4.

After writing this paper, we came across an older, nonperiodic example due to Marjorie Rice [R], pictured in Figure 5.

Meanwhile, by trial and error with Geometer's Sketchpad, we found other such Cario-Prismatic tilings with symmetries of four of the seventeen wallpaper groups (Figures 6−10) and others with fewer or no symmetries (Figures 11−15) ); see [C] for more. We wonder whether there are examples for the other wallpaper groups. For a chart of the 17 wallpaper groups and their respective symmetries, refer to [S2] or [W].

Daniel Huson [Hu] has found tilings by many other pairs of prototiles.

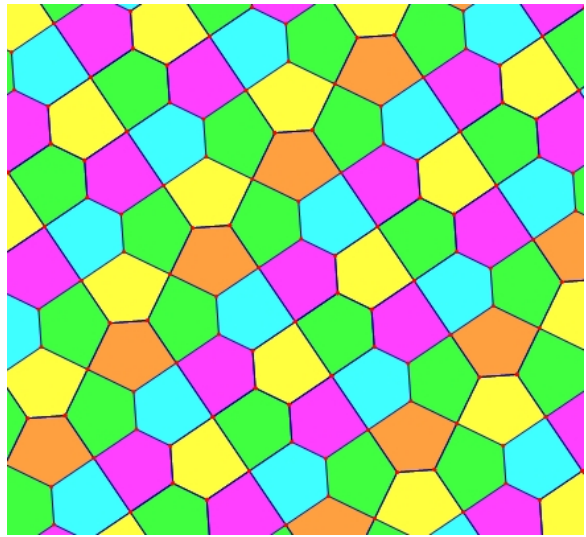

Figure 4. Another tiling by a mixture of Cairo and Prismatic tiles, in which the diagonal rows of Cairo tiles are separated by three diagonal rows of Prismatic tiles.

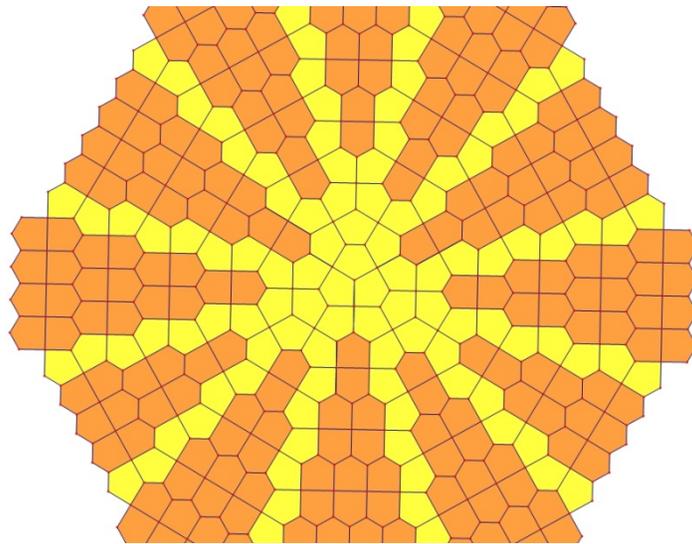

Figure 5. A nonperiodic tiling by a mixture of Cairo and Prismatic tiles due to Marjorie Rice [R]

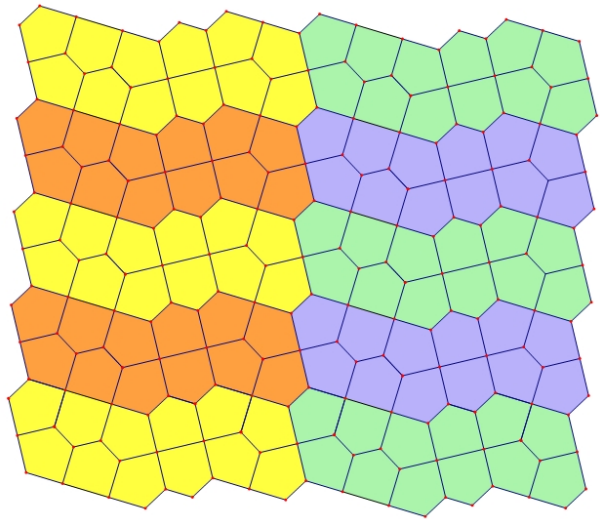

Figure 6. Cairo-Prismatic tiling with wallpaper group p1

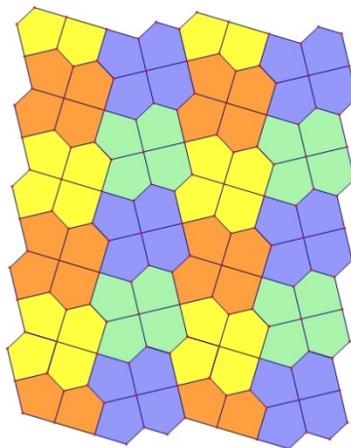

Figure 7. Cairo-Prismatic tiling with wallpaper group p2

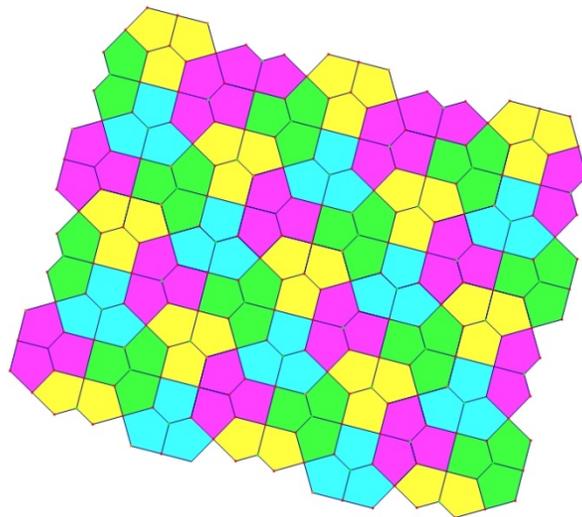

Figure 8. Cairo-Prismatic tiling with wallpaper group p4g

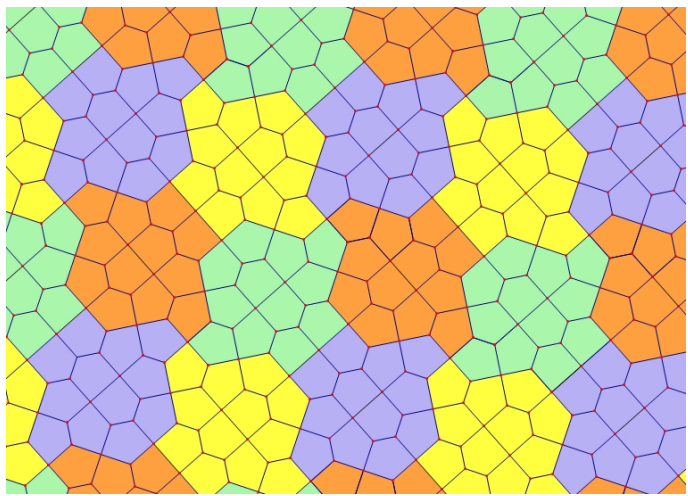

Figure 9. Another Cairo-Prismatic tiling with wallpaper group p4g

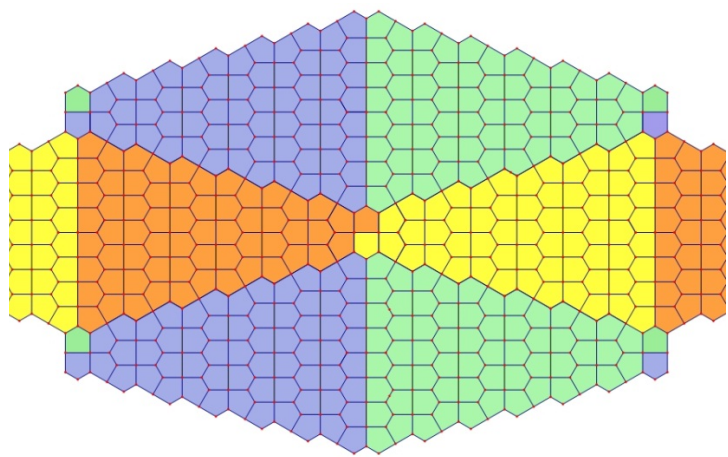

Figure 10. Cairo-Prismatic tiling with wallpaper group cmm

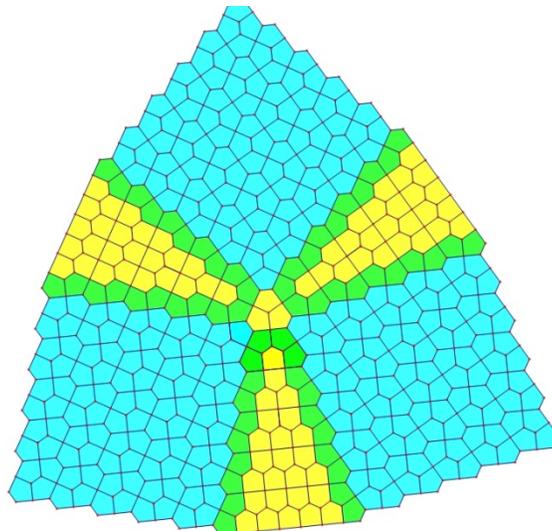

Figure 11. "Windmill" Cairo-Prismatic tiling with symmetry group D₃

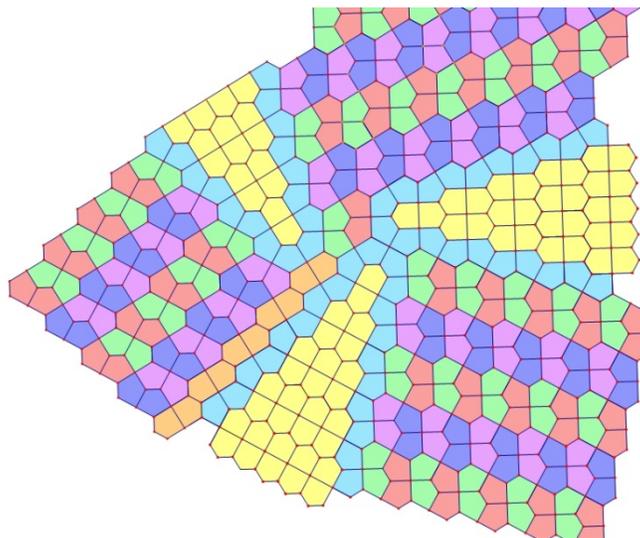

Figure 12. "Chaos" Cairo-Prismatic tiling with no symmetry

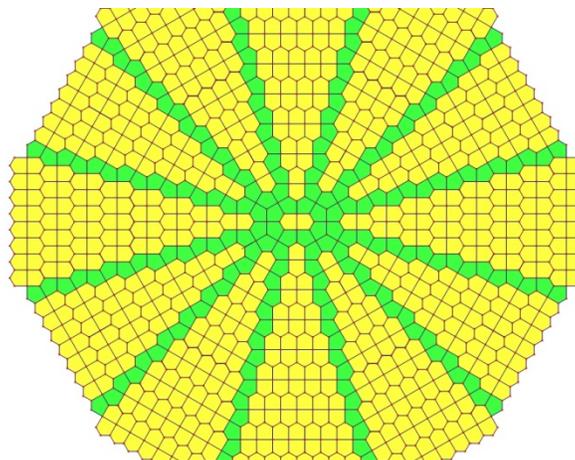

Figure 13. "Plaza" Cairo-Prismatic tiling with symmetry group D₂

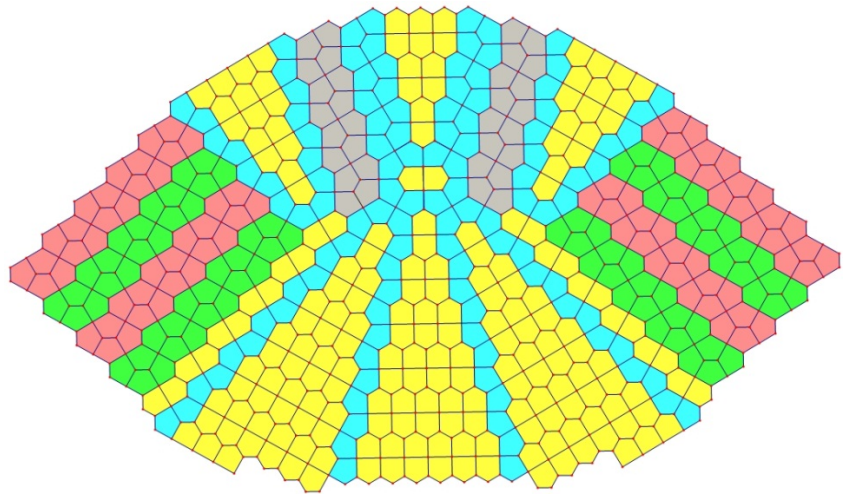

Figure 14. "Bunny" Cairo-Prismatic tiling with symmetry group D₁

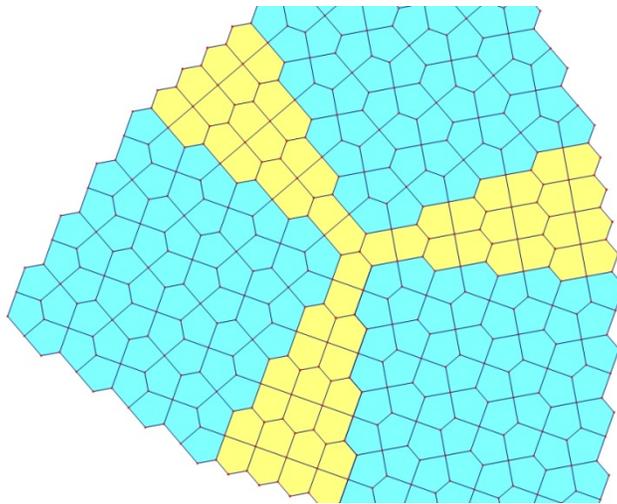

Figure 15. "Waterwheel" Cairo-Prismatic tiling with symmetry group D₃

# 3. Isoperimetric Pentagonal Tilings

Our main Theorem 3.5 proves that the Cairo and Prismatic tilings minimize perimeter among all tilings by unit-area convex polygons with at most five sides. We begin with the isoperimetric Proposition 3.1 which we have not been able to find in the literature.

**Proposition 3.1.** *For n given angles $0 < a_i \leq \pi$ summing to $(n-2)\pi$, the n-gon circumscribed about the unit circle as in Figure 16 is uniquely perimeter minimizing for its area. Scaled to unit area, its perimeter is*

*(1)*
$$2\sqrt{\sum \cot(a_i/2)}\,.$$

*Since cotangent is strictly convex up to $\pi/2$, the more nearly equal the angles, the smaller the perimeter.*

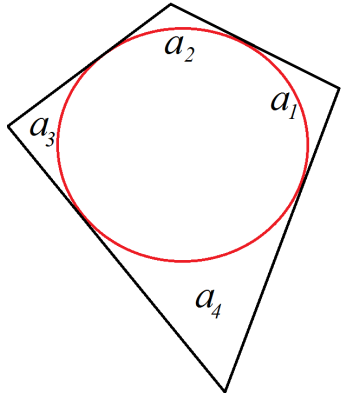

Figure 16. For given angles the least-perimeter unit-area polygon is circumscribed about a circle.

*Proof.* For any (not necessarily even) norm on $\mathbf{R}^2$ (or $\mathbf{R}^n$), the unit ball in the dual norm or "Wulff shape" uniquely minimizes the integral of the norm of the unit normal over the boundary for given area (see [M2, 10.6]). In particular, the circumscribed *n*-gon minimizes the integral of a norm which is 1 on its normals and some irrelevant (larger) value on other directions. By simple geometry, its perimeter $P_0$ and area $A_0$ satisfy

$$P_0 = 2A_0 = 2\sum \cot(a_i/2)\,.$$

The final formula for unit area follows by scaling.

*Remarks.* If you allow angles between $\pi$ and $2\pi$, the circumscribed polygon still minimizes perimeter for sides in the associated directions, but they may necessarily occur in a different order and produce different angles and a smaller perimeter. If you require the prescribed order to keep the prescribed angles, some sides may become infinitesimally small, thus omitted, replacing pairs of external angles with their sum and reducing our equation (1). For n > 3, the least-perimeter unit-area *n*-gon with one small prescribed angle is not convex; for n = 4 it is close to being an equilateral triangle with a small nonconvex protrusion near one of the three vertices, as in Figure 17.

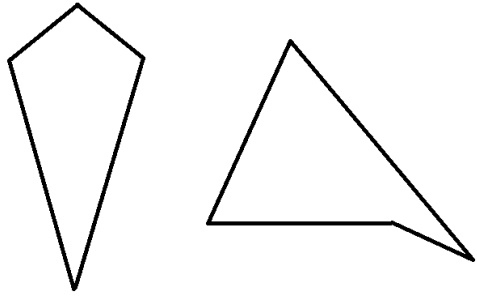

Figure 17. The least-perimeter unit-area quadrilateral with one small prescribed angle is not the convex quadrilateral with the other angles equal but a triangle with a small nonconvex protrusion.

**3.2. Lemma.** *For positive integers k + k' = n, consider the positive function*

$$f^2(u) = g(u) = kt + k't'$$

*where t = cot θ, t' = cot θ', 0 < θ, θ' < π/2, θ = k'u, θ' = (n-2)π/2k' - ku. Then f is strictly convex.*

*Proof.* It suffices to show that $2gg'' - g'^2$ is positive (where primes on *g* denote differentiation with respect to u). Writing *s = csc θ, s' = csc θ'*, we compute that

$$g' = -kk's^2 + kk's'^2,$$

$$g'' = 2kk'^2s^2t + 2k^2k's'^2t',$$

and after discarding some positive terms

$$2gg'' - g'^2 > 2k^2k'^2s^2t^2 + 2k^2k'^2s'^2t'^2 - k^2k'^2s^4 - k^2k'^2s'^4 + 2k^2k'^2s^2s'^2,$$

which is proportional to

$$2s^2t^2 + 2s'^2t'^2 - s^4 - s'^4 + s^2s'^2$$

$$> 2s^2t^2 + 2s'^2t'^2 - s^4 - s'^4 + s^2 + s'^2$$

$$= 2s^2t^2 + 2s'^2t'^2 - s^2t^2 - s'^2t'^2 > 0.$$

**3.3. Lemma.** *For positive constants n, c, c' and nonnegative constants d, d', consider the positive function*

$$f^2(k) = g(k) = kt + k't',$$

*where $0 < k < n$, $k' = n-k$, $t = \cot \theta$, $t' = \cot \theta'$, $0 < \theta, \theta' < \pi/2$, $\theta = d - c/k$, $\theta' = d - c'/k'$, and $\min\{\theta, \theta'\} \le \pi/3$. Then $f$ is strictly convex.*

*Proof.* It suffices to show that $2gg''-g'^2$ is positive. Write $s = \csc \theta$, $s' = \csc \theta'$. By symmetry we may assume that $\theta \le \theta'$, so $c/k \ge c'/k'$, $t \ge t'$, and $s \ge s'$. We compute that

$$g' = t - ks^2(c/k^2) - t' + k's'^2(c'/k'^2) = t - t' - (s^2c/k - s'^2c'/k'),$$

so

$$-s^2c/k \le g' \le t - t'.$$

Since

$$t - t' \le s^2(\theta-\theta) \le s^2c/k,$$

we see that

$$|g'| \le s^2c/k.$$

We further compute that

$$g'' = 2s^2tc^2/k^3 + 2s'^2tc'^2/k^3.$$

Hence

$$2gg'' - g'^2 > 2(kt)(\, 2s^2tc^2/k^3) - s^4c^2/k^2,$$

which is proportional to $4t^2 - s^2$, which is nonnegative because $\theta \le \pi/3$.

**Corollary 3.4.** *For integers $n \ge 3$, the positive function $f_n$ defined by*

$$f_n^2(k, k\theta) = k \cot \theta + k' \cot \theta',$$

*where $k'$ and $\theta'$ are defined by $k + k' = n$, $k\theta + k'\theta' = 3\pi/2$, $0 < k, k' < n$, $0 < \theta, \theta' < \pi/2$, and $\min\{\theta, \theta'\} \le \pi/3$ is strictly convex.*

*Proof.* On a vertical line in the domain ($k$ constant), $f_n$ is strictly convex by Lemma 3.2. On a nonvertical line, $f_n$ is strictly convex by Lemma 3.3.

**Theorem 3.5.** *Perimeter-minimizing planar tilings by unit-area convex polygons with at most five sides are given by Cairo and Prismatic tiles (as in Figures 1, 2, 4-15).*

*Remarks.* There is never uniqueness, since compact variations do not change the limiting perimeter ratio.

For doubly periodic or monohedral tilings, however, there are no others, although the tilings need not be edge to edge.

If nonconvex tiles are allowed, the result remains conjectural (except for monohedral tilings, when such a large angle gives the prototile more perimeter than a square).

We are not assuming that the tiling is edge to edge.

*Proof of Theorem 3.5.* Consider a planar tiling with bounded perimeter ratio ρ by unit-area convex polygons with at most five sides. We may assume that every angle is strictly less than π by eliminating vertices at angles of π from polygons. (We are not assuming that the tilings are edge to edge.)

To control truncation, we claim that for some sequence of $R$ going to infinity, the circle about the origin of radius $R$ meets the tiling in o($R^2$) points (*i.e.*, the ratio to $R^2$ goes to 0). Otherwise for almost all large $R$, for some ε > 0, since $P$ grows at a rate at least equal to the number of such points (see [M1, 15.3] for details),

$$dP/dr \geq \varepsilon R^2 \geq (\varepsilon/\rho)P,$$

which implies that perimeter and hence area grow exponentially, a contradiction, since area equals $\pi R^2$.

For such a sequence of R, consider the tiles inside the disc $D_R$. The ratio of their area to $\pi R^2$ approaches 1, while the number of adjacent discarded tiles is o($R^2$), negligible. Therefore the perimeter ratio inside $D_R$ equals the average cost of the tiles inside $D_R$ up to o(1).

We will examine more carefully our sequence as $R$ approaches infinity. By taking a subsequence if necessary, we may assume that the fraction of polygons of $n$ edges approaches a limit $f_n$, hence that the number of edges per polygon approaches $e = \Sigma n f_n$, and that the number of vertices per polygon of degree $d$ in the tiling approaches a limit $v_d$. Thus

(1)      $\Sigma v_d = e.$

By Euler, for the tiles inside $D_R$,

(2)      $V = E - F + 1 \geq eF/2 - F + 1 = (e/2 - 1)F + 1,$

with equality up to o(1) if the tiling is edge to edge. Hence, if the tiling is edge to edge, the following equality holds up to o(1) for the tiles inside $D_R$:

(3)      $\Sigma v_d / d = V/F = e/2 - 1.$

If the tiling is not edge to edge, we count a vertex where $k-1$ other edges meet in the interior of an edge as contributing to $v_{2k}$ instead of $v_k$, because these polygon angles are $2\pi/2k$ on the average. Such a vertex counts just 1/2 on the left side of (3). Let $V^*$ be the number of such vertices. Then

(4)      $\Sigma v_d / d = (V - V^*/2)/F = (E - V^*/2 - F)/F = ((eF + V^*)/2 - V^*/2 - F)/F = e/2 - 1,$

because each such vertex breaks a polygonal edge into two edges, replacing $eF$ with $eF + V^*$. In the limit, equation (4) holds exactly. By (1) and (4),

$$v_3/3 + (e-v_3)/4 \geq e/2 - 1,$$

(5)     $v_3 \geq 3e - 12.$

Now it suffices to prove for every collection of $N$ convex unit-area polygons with at most five edges with $(3e-12)N - N_1$ angles of average measure at least $2\pi/3$ with $N_1 =$ o($N$), the average perimeter is at least nearly that of the Cairo and Prismatic pentagons. Among all such collections of $N$ pentagons, choose one to minimize the average perimeter. By Proposition 3.1, they are all circumscribed about circles. Since by Proposition 3.1 reducing large angles reduces perimeter, the $(3e-12)N - N_1$ "large" angles have average measure exactly $2\pi/3$. Call the remaining angles the "small" angles, all at most $2\pi/3$. Up to o($N$)$N$, the number of small angles is $(12-2e)N$, their sum is $(6-e)\pi N$, so their average is $\pi/2 + $ o($N$) and all large angles are at least that large.

For each $3 \leq n \leq 5$ we may assume that the average large angle in an $n$-gon approaches a limit $2\theta_n$ and that the number of large angles per $n$-gon approaches a limit $k_n$. Then

(6)     $$\Sigma f_n k_n = 3e-12, \ \Sigma f_n k_n \theta_n = (\pi/3) \ \Sigma f_n k_n.$$

Since $e = \Sigma n f_n$,

(7)     $f_3(k_3+3) + f_4 k_4 + f_5(k_5-3) = 0, f_3(k_3\theta_3+\pi) + f_4 k_4\theta_4 + f_5(k_5\theta_5-\pi) = 0.$

Note that each quadrilateral has at least one large angle, or it could be replaced by a Cairo pentagon, increasing both $(3e-12)N - N_1$ and the number of large angles by three and decreasing perimeter.

By Proposition 3.1, in each $n$-gon, the large angles are all equal and the small angles are all equal.

By convexity (Corollary 3.4), for each $3 \leq n \leq 5$, in the limit the average perimeter is bounded below by the values of the perimeter function $P_n(k,q)$ at $k_n$ and $k_n\theta_n$, where

$$P_n(k, q) = 2(k \cot q/k + k' \cot q'/k')^{1/2},$$

$k + k' = n, \ q + q' = (n-2)\pi/2.$ We compute to four decimal places that

(8)     $(\partial P_3/\partial k)(1/2, \pi/8) \approx -.7217$ ,     $(\partial P_3/\partial q)(1/2, \pi/8) \approx 1.2265$ ,

   $(\partial P_4/\partial k)(1, \pi/3) \approx -.4455$ ,     $(\partial P_4/\partial q)(1, \pi/3) \approx .5334$ ,

   $(\partial P_5/\partial k)(3, \pi ) \approx -.3091$ ,     $(\partial P_5/\partial q)(3, \pi) \approx .3451$ .

Note that $f_5 > 0$; otherwise by (7) $f_4 = 1$, $k_4 = 0$, and the square tiling provides a lower bound, while we know that the Cairo tiling is better. Also by (7) $k_5 \leq 3$. Since the average small angle is $\pi/2$, every large angle $2\theta_n$ is at least $\pi/2$, $k_3 < 2$, and $k_4 < 3$. Since the sum of the angles of an $n$-gon is $(n-2)\pi$, $k_n\theta_n + k_n'\theta_n' = (n-2)\pi/2$, $\theta_5 \geq 3\pi/10$ and $\theta_n' \leq (n-2)\pi/2n$.

Since by Corollary 3.4 each $f_n$ is convex, we can bound the perimeter $P$ per tile from below by the linear approximations:

$P = f_3 P_3(k_3, k_3\theta_3) + f_4 P_4(k_4, k_4\theta_4) + f_5 P_5(k_5, k_5\theta_5)$

$\geq\ f_3 \left[P_3(1/2, \pi/8) + (\partial P_3/\partial k)(1/2, \pi/8)(k_3-1/2) + (\partial P_3/\partial q)(1/2, \pi/8)(k_3\theta_3-\pi/8)\right]$
$+ f_4 \left[P_4(1, \pi/3) + (\partial P_4/\partial k)(1, \pi/3)(k_4-1) + (\partial P_4/\partial q)(1, \pi/3)(k_4\theta_4-\pi/3)\right]$
$+ f_5 \left[P_5(3, \pi) + (\partial P_5/\partial k)(3, \pi)(k_5-3) + (\partial P_5/\partial q)(3, \pi)(k_5\theta_5-\pi)\right].$

Using (7) to substitute for $f_5(k_5-3)$ and $f_5(k_5\theta_5-\pi)$ we obtain

$P \geq f_5 P_5(3, \pi) + f_3 Q_3 + f_4 Q_4\,,$

where (using appropriate roundings up or down)

$Q_3 > P_3(1/2, \pi/8) + (1/2)(.7216)\ - (\pi/8)(1.2266) + 3(.3090) - \pi(.3452)$
$+ k_3(-.7218 + .3090) + k_3\theta_3(1.2264 - .3452)$

$> 4.6503 + (1/2)(.7216)\ - (\pi/8)(1.2266) + 3(.3090) - \pi(.3452)$
$+ k_3(-.7218 + .3090) + k_3\theta_3(1.2264 - .3452)$

$> 4.3718 + k_3(-.4128 + (\pi/4)(.8812)) > P_5(3, \pi)$

and

$Q_4 > P_4(1, \pi/3) + .4454\ - (\pi/3)(.5335) + k_4(-.4456 + .3090) + k_4\theta_4(.5333 - .3452)$

$> 3.9622 + k_4(-.1366 + (\pi/4)(.1881)) > P_5(3, \pi).$

Therefore the perimeter per tile is at least the perimeter $P_5(3, \pi)$ of the Cairo and the Prismatic tile, as desired.

For doubly periodic tilings, the argument can be done on the torus without limits or truncation error, and every tile is Cairo or Prismatic.

**Remark 3.6.** Theorem 3.5 holds for a monohedral edge-to-edge planar tiling by a *curvilinear* polygon $\mathscr{P}$ (embedded rectifiable curves meeting only at vertices) with at most five edges. Because the tiling is monohedral and edge-to-edge, one can replace $\mathscr{P}$ with an immersed rectilinear $n$-gon with the same vertices. If it is not embedded, it has more perimeter than a square.

Ping Ngai Chung
Department of Mathematics
MIT
briancpn@mit.edu

Miguel A. Fernandez
Mathematics and Computer Science Department
Truman State University
maf2831@truman.edu

Yifei Li
Department of Mathematics and Computer Science
Berea College
Yifei1124@gmail.com

Michael Mara
Department of Mathematics and Statistics
Williams College
mtm1@williams.edu

Frank Morgan
Department of Mathematics and Statistics
Williams College
Frank.Morgan@williams.edu

Isamar Rosa
Department of Mathematical Sciences
University of Puerto Rico at Mayaguez
silver.rose.isa@gmail.com

Niralee Shah
Department of Mathematics and Statistics
Williams College
nks1@williams.edu

Luis Sordo Vieira
Department of Mathematics
Wayne State University
dw8603@wayne.edu

Elena Wikner
Department of Mathematics and Statistics
Williams College
elena.wikner@gmail.com